\documentclass[12pt]{article}
\usepackage{amsfonts,amsmath,amsxtra}
\usepackage{amssymb,amscd}
\usepackage[mathscr]{eucal}

\def\hybrid{\topmargin 0pt      \oddsidemargin 0pt
        \headheight 0pt \headsep 0pt
        \textwidth 16.5cm
        \textheight 23cm
        \marginparwidth 0.0in
        \parskip 5pt plus 1pt   \jot = 1.5ex}
\catcode`\@=11
\def\marginnote#1{}
\newcount\hour
\newcount\minute
\newtoks\amorpm
\hour=\time\divide\hour by60
\minute=\time{\multiply\hour by60 \global\advance\minute by-\hour}
\edef\standardtime{{\ifnum\hour<12 \global\amorpm={am}%
        \else\global\amorpm={pm}\advance\hour by-12 \fi
        \ifnum\hour=0 \hour=12 \fi
      \number\hour:\ifnum\minute<10 0\fi\number\minute\the\amorpm}}
\edef\militarytime{\number\hour:\ifnum\minute<10 0\fi\number\minute}

\def\draftlabel#1{{\@bsphack\if@filesw {\let\thepage\relax
   \xdef\@gtempa{\write\@auxout{\string
      \newlabel{#1}{{\@currentlabel}{\thepage}}}}}\@gtempa
   \if@nobreak \ifvmode\nobreak\fi\fi\fi\@esphack}
        \gdef\@eqnlabel{#1}}
\def\@eqnlabel{}
\def\@vacuum{}
\def\draftmarginnote#1{\marginpar{\raggedright\scriptsize\tt#1}}

\def\draft{\oddsidemargin -0.1truein
        \def\@oddfoot{\sl preliminary draft \hfil
        \rm\thepage\hfil\sl\today\quad\militarytime}
        \let\@evenfoot\@oddfoot \overfullrule 3pt
        \let\label=\draftlabel
        \let\marginnote=\draftmarginnote
\def\@eqnnum{{\rm (\theequation)}
\rlap{\kern\marginparsep\tt\@eqnlabel}%
\global\let\@eqnlabel\@vacuum}  }

\newfont{\Bbbb}{msbm7 scaled 1\@ptsize00}
\newcommand{\zs}{\raise-1pt\hbox{$\mbox{\Bbbb Z}$}}

\font\sevenmsa=msam6 
\newfam\msafam
\textfont\msafam=\sevenmsa
\def\hexnumber@#1{\ifnum#1<10 \number#1\else
\ifnum#1=10 A\else\ifnum#1=11 B\else\ifnum#1=12 C\else
\ifnum#1=13 D\else\ifnum#1=14 E\else\ifnum#1=15 F\fi\fi\fi\fi\fi\fi\fi}
\def\msa@{\hexnumber@\msafam}
\def\llcorner{\delimiter"4\msa@78\msa@78 }
\def\lrcorner{\delimiter"5\msa@79\msa@79 }
\mathchardef\blacktriangleright="3\msa@49
\mathchardef\blacktriangleleft="3\msa@4A
\font\tenmsb=msbm10 scaled 1\@ptsize00
\newfam\msbfam
\textfont\msbfam=\tenmsb
\scriptfont\msbfam=\tenmsb

\newdimen\linethick  \linethick=0.4pt
\newdimen\hboxitspace    \hboxitspace=5pt
\newdimen\vboxitspace    \vboxitspace=5pt

\def\fr#1{%
\be\new
\vcenter{
\hrule height\linethick
           \hbox{\vrule width\linethick
                 \kern\hboxitspace
                 \vbox{\kern\vboxitspace
                       \hbox{$\begin{array}{c}\displaystyle#1
          \end{array}$}%
                       \kern\vboxitspace}%
                 \kern\hboxitspace
                 \vrule width\linethick}%
           \hrule height\linethick}%
\ee}

\newdimen\Squaresize \Squaresize=14pt
\newdimen\Thickness \Thickness=0.5pt

\def\Square#1{\hbox{\vrule width \Thickness
   \vbox to \Squaresize{\hrule height \Thickness\vss
      \hbox to \Squaresize{\hss#1\hss}
   \vss\hrule height\Thickness}
\unskip\vrule width \Thickness}
\kern-\Thickness}

\def\Vsquare#1{\vbox{\Square{$#1$}}\kern-\Thickness}

\def\numberbysection{\@addtoreset{equation}{section}
        \def\theequation{\thesection.\arabic{equation}}}
\numberbysection

\renewcommand{\theequation}{\thesection.\arabic{equation}}
\def\titlepage{\@restonecolfalse\if@twocolumn\@restonecoltrue\onecolumn
     \else \newpage \fi \thispagestyle{empty}\c@page\z@
        \def\thefootnote{\fnsymbol{footnote}} }

\def\endtitlepage{\if@restonecol\twocolumn \else  \fi
        \def\thefootnote{\arabic{footnote}}
        \setcounter{footnote}{0}}  
\relax

\hybrid
\makeatletter
\newdimen\normalarrayskip            
\newdimen\minarrayskip               
\normalarrayskip\baselineskip
\minarrayskip\jot
\newif\ifold             \oldtrue            \def\new{\oldfalse}
\def\arraymode{\ifold\relax\else\displaystyle\fi}
\def\eqnumphantom{\phantom{(\theequation)}} 
\def\@arrayskip{\ifold\baselineskip\z@\lineskip\z@
     \else
     \baselineskip\minarrayskip\lineskip1\baselineskip\fi}

\def\@arrayclassz{\ifcase \@lastchclass \@acolampacol \or
\@ampacol \or \or \or \@addamp \or
   \@acolampacol \or \@firstampfalse \@acol \fi
\edef\@preamble{\@preamble
  \ifcase \@chnum
     \hfil$\relax\arraymode\@sharp$\hfil
     \or $\relax\arraymode\@sharp$\hfil
     \or \hfil$\relax\arraymode\@sharp$\fi}}

\def\@array[#1]#2{\setbox\@arstrutbox=\hbox{\vrule
     height\arraystretch \ht\strutbox
     depth\arraystretch \dp\strutbox
width\z@}\@mkpream{#2}\edef\@preamble{\halign \noexpand\@halignto
\bgroup \tabskip\z@ \@arstrut \@preamble \tabskip\z@ \cr}%
\let\@startpbox\@@startpbox \let\@endpbox\@@endpbox
  \if #1t\vtop \else \if#1b\vbox \else \vcenter \fi\fi
  \bgroup \let\par\relax
  \let\@sharp##\let\protect\relax
  \@arrayskip\@preamble}
\def\eqnarray{\stepcounter{equation}%
              \let\@currentlabel=\theequation
              \global\@eqnswtrue
              \global\@eqcnt\z@
              \tabskip\@centering              
              \let\\=\@eqncr
              $$%
            \halign to \displaywidth  \bgroup
             \eqnumphantom \@eqnsel
      \hskip\@centering                               
    $\displaystyle  \tabskip\z@ {##}$%
    &\global\@eqcnt\@ne \hskip 2\arraycolsep
         $ \displaystyle  \arraymode{##}$\hfil
    &\global\@eqcnt\tw@ \hskip 2\arraycolsep
         $\displaystyle\tabskip\z@{##}$\hfil
         \tabskip\@centering
    &{##}\tabskip\z@\cr}
\makeatother
\newtheorem{te}{Theorem}[section]

\newtheorem{prop}{Proposition}[section]           

\newtheorem{ex}{Example}[section]

\newcommand{\beq}[1]{\begin{equation}\label{#1}}
\newcommand\eeq{\end{equation}}
\newcommand\bqa{\begin{eqnarray}}
\newcommand\eqa{\end{eqnarray}}
\def\be{\begin{eqnarray}\new\begin{array}{cc}}
\def\ee{\end{array}\end{eqnarray}}

\def\beq{\begin{equation}}
\def\eeq{\end{equation}}
\def\bse{\begin{subequations}}                
\def\ese{\end{subequations}}
\def\bp{\begin{pmatrix}}
\def\ep{\end{pmatrix}}

\def\i{\imath}

\def\stack#1#2{\raise0.7pt\hbox{$\mathrel{\mathop{#2}\limits^{#1}}$}}
\def\tr{\triangleright}
\def\tl{\triangleleft}
\def\sem{\mathsurround=0pt \raise1pt
\hbox{$\scriptscriptstyle>\!\!$}\:\!\!\tl}
\def\mes{\mathsurround=0pt \tr\!\:\!\raise0.8pt
\hbox{$\scriptscriptstyle\!\!<$}\,}
\def\]{\mathsurround=0pt ]\raise-2pt\hbox{$_\ast$}}

\def\la{\lambda}

\def\<{\langle}
\def\>{\rangle}

\def\cl{{\cal L}}

\def\we{\raise-1pt\hbox{$\,\stackrel{\wedge}{,}\,$}}

\def\pr {\partial}

0
\begin{document}

\footnotesize
\normalsize

\newpage

\thispagestyle{empty}

\begin{center}

\phantom.
\bigskip
{\hfill{\normalsize hep-th/0609082}\\
\hfill{\normalsize HMI-06-03}\\
\hfill{\normalsize ITEP-TH-06-42}\\
\hfill{\normalsize TCD-MATH-06-XX}\\
[10mm]\Large\bf
Baxter $Q$-operator and Givental integral representation  for
 $C_n$ and $D_n$}
\vspace{0.5cm}

\bigskip\bigskip
{\large A. Gerasimov}
\\ \bigskip
{\it Institute for Theoretical and
Experimental Physics, Moscow, 117259, Moscow,  Russia} \\
{\it  School of Mathematics, Trinity
College, Dublin 2, Ireland } \\
{\it  Hamilton
Mathematics Institute, TCD, Dublin 2, Ireland},\\
\bigskip
{\large D. Lebedev\footnote{E-mail: lebedev@ihes.fr}}
\\ \bigskip
{\it Institute for Theoretical and Experimental Physics, 117259,
Moscow, Russia}  \\
{\it l'Institute des Hautes
\'Etudes Scientifiques, 35 route de Chartres,
Bures-sur-Yvette, France},\\
\bigskip
{\large S. Oblezin} \footnote{E-mail: Sergey.Oblezin@itep.ru}\\
\bigskip {\it Institute for Theoretical and Experimental Physics,
117259, Moscow,
Russia}\\
{\it Max-Planck-Institut f\"ur Mathematik, Vivatsgasse 7, D-53111
Bonn, Germany},\\
\end{center}

\vspace{0.5cm}

\begin{abstract}
\noindent

Recently  integral representations for  the 
 eigenfunctions of  quadratic open Toda chain Hamiltonians
 for classical groups was proposed. This representation 
generalizes Givental representation for $A_n$.
In this note we verify that the wave functions defined by these 
 integral representations are  common eigenfunctions for the 
complete  set of open Toda chain Hamiltonians. We consider the zero eigenvalue 
wave functions for  classical groups $C_n$ and $D_n$ thus completing 
the generalization of the Givental construction in these cases.    
The construction is based on a recursive procedure
and uses the formalism of Baxter $Q$-operators. We also verify 
 that the  integral $Q$-operators for $C_n$, $D_n$ and twisted affine algebra
$A_{2n-1}^{(2)}$ proposed previously intertwine  complete  sets of  
Hamiltonian operators.  Finally we provide  integral
  representations of the eigenfunctions of  the 
quadratic $D_n$ Toda chain Hamiltonians for generic nonzero eigenvalues.

\end{abstract}

\vspace{1cm}

\clearpage \newpage


\normalsize
\section{Introduction}

In \cite{GLO}   we provided  integral representations for
 zero eigenvalue wave functions of open Toda chain quadratic Hamiltonian
operators corresponding to the classical series of finite Lie algebras $B_n$, $C_n$,
and $D_n$. This representation  generalizes Givental representation for $A_n$ \cite{Gi}.
  In this short note we  verify the eigenvalue properties for the full set of the
Hamiltonian operators for $C_n$ and $D_n$ open  Toda chains.
 We also give  as an  example of the  
integral representations for generic eigenvalues 
the explicit expression for the wave function of $D_n$ Toda chain.

We use a recursive representation of the 
 common eigenfunctions of the full sets of $C_n$ and $D_n$ Toda chain
Hamiltonians. The integral  recursive operator
is naturally represented as a composition of two  elementary integral operators
with simple integral kernels. The interesting property of 
these elementary integral operators is that they relate $C_n$ and
$D_n$  wave functions. This is the main reason why we consider these two cases together.

Baxter $Q$-operator is a key tool in  the theory of quantum
integrable systems \cite{B}. It plays the role of the generator of quantum
Backlund transformations and was explicitly constructed  for $A_{n}^{(1)}$
closed Toda chain in \cite{PG}.  It was noted in 
 \cite{GKLO}  that the recursive operators for $A_{n}$ open Toda chain 
 obtained as a limit of   $Q$-operators of $A_{n}^{(1)}$ Toda chain
 play the crucial role in the construction of the 
Givental integral representation.
 In \cite{GLO} we generalize this relation by constructing
$A_{2n-1}^{(2)}$ Baxter $Q$-operator and demonstrating  that $C_n$
and $D_n$ recursive operators can be obtained from $A_{2n-1}^{(2)}$
Baxter $Q$-operator in a similar limit. The properties of the
recursive operators are then follows from the commutation
relation between quantum $L$-operator and $Q$-operator
of  $A^{(2)}_{2n-1}$ closed Toda chain.
In this note we verify that the integral $Q$-operators for 
$A_{2n-1}^{(2)}$ Toda chain introduced in \cite{GLO} commute with 
the full set of $A_{2n-1}^{(2)}$ Toda chain Hamiltonians. 
The proof is based on the explicit commutation relations with 
$L$-operator for Toda chain. In a certain limit this leads 
to the commutation relations between recursive operators 
 and Hamiltonian operators of $C_n$ and $D_n$ open Toda chains. 
The  fact that  the wave functions constructed in \cite{GLO} 
are common eigenfunctions of the full set of the Toda chain 
then easily follows form the commutation relations. 
The use of $A_{2n-1}^{(2)}$ is not very essential in this construction
and is dictated by the simplicity of the resulting degeneration
procedure. For example similar
relations exist between $Q$-operators for  $C^{(1)}_n$ and
$D_n^{(1)}$   and recursive operators for $C_n$ and $D_n$. These
constructions will be presented  elsewhere together with a complete account of the
constructions of the Baxter $Q$-operators, recursive operators and
integral representations of wave functions for all classical
series \cite{GLO1}.

The plan of this note  is as follows. In Section 2
we  recall  explicit expressions for  the Baxter $Q$-operator
for $A_{2n-1}^{(2)}$ given in \cite{GLO} and prove that it commutes with
the full set of  $A_{2n-1}^{(2)}$ Toda chain Hamiltonian operators.
In  Section 3 we derive recursive operators
 for $C_n$ and $D_n$ as a certain limit of the $A_{2n-1}^{(2)}$
 Baxter $Q$-operator. We prove the intertwining relations 
of the recursive operators with Hamiltonians of 
$C_n$ and $D_n$ Toda chains.  In Section 4
integral representations of  the 
eigenfunctions  of $D_n$  Toda chain quadratic Hamiltonian
operators for generic eigenvalues are given. 

Finally in the Appendix  we explicitly check  that thus obtained
wave functions for  $D_2$ open Toda chain can be reduced
to  a product of the wave function of two independent $A_1$
Toda chains. Although it is a consequence of the isomorphism
 $D_2=A_1\oplus A_1$ the transformation is rather nontrivial
and serves as an independent check of our construction.

{\em Acknowledgments}:    The research of A.~Gerasimov
is partly supported by the grant RFBR 04-01-00646
and Enterprise Ireland Basic Research Grant;
The research of D.~Lebedev is partially supported
by the grant RFBR 04-01-00646.
The research of S. Oblezin was partially supported by the grant RFBR
04-01-00642.   D.~Lebedev is grateful to
 Institute des Hautes \'Etudes Scientifiques for a warm
hospitality. S.~Oblezin is grateful to Max-Planck-Institut f\"ur Mathematik
for excellent working conditions.

\section{Baxter $Q$-operator for  $A_{2n-1}^{(2)}$}

In this section we recall the construction of  $Q$-operator
for $A_{2n-1}^{(2)}$ and explicitly check that it commutes with
the full set of $A_{2n-1}^{(2)}$ closed Toda chain Hamiltonians.
Let us start with the description of  $A_{2n-1}^{(2)}$ closed
Toda chain. For the relevant facts on the root systems see \cite{K},\cite{DS}. The
detailed description of the  corresponding Toda chains can be found
in \cite{RSTS}.

Let us fix an orthonormal  basis $\{e_1,\ldots,e_{n}\}$ in
$\mathbb{R}^{n}$. Simple roots of the twisted affine root system $A_{2n-1}^{(2)}$
can be the represented in the following form
 \be \alpha_1=2e_1,\qquad
\alpha_{i+1}=e_{i+1}-e_i,\qquad 1\leq i\leq n-1\qquad
\alpha_{n+1}=-e_n-e_{n-1}, \ee
The corresponding Dynkin diagram is given by
$$
\begin{CD}\alpha_1 \Longrightarrow  \alpha_2 @>>> \ldots @>>> \alpha_{n-1}
@>>> \alpha_n\\ @. @. @VVV @.\\ @. @.\alpha_{n+1}@.\end{CD}
$$

To construct $L$-operator of the closed Toda chain
corresponding to an affine root system  one should chose an evaluation
representation of the corresponding affine Lie algebra.
Let $\{e_i,h_i,f_i\}$ be  a bases of the twisted Lie
algebra  $\mathfrak{g}$ corresponding to root system $A_{2n-1}^{(2)}$.
We define the evaluation representation of $\widehat{\mathfrak{g}}$ as follows
(see e.g. \cite{DS}). Consider  the following  homomorphism
\be
\pi:\widehat{\mathfrak{g}}\longrightarrow Mat(2n,\mathbb{C})\otimes \mathbb{C}[u,u^{-1}]
\ee
defined  explicitly as
 \be
\pi(e_0)=u(E_{1,2n-1}+E_{2,2n})/2,\qquad \pi(f_0)=2u^{-1}(E_{2n-1,1}+E_{2n,2}),\\
\pi(e_i)=E_{i+1,i}+E_{2n+1-i,2n-i},\qquad \pi(f_i)=E_{i,i+1}+E_{2n-i,2n+1-i},\\
\pi(e_n)=E_{n+1,n},\qquad \pi(f_n)=2E_{n,n+1},\\
\pi(h_0)=(E_{1,1}-E_{2n,2n})+(E_{2,2}-E_{2n-1,2n-1}),\\
\pi(h_i)=(E_{i+1,i+1}-E_{i,i})+(E_{2n+1-i,2n+1-i}-E_{2n-i,2n-i}),\\
 \pi(h_n)=2(E_{n+1,n+1}-E_{n,n})\ee where
$(E_{ij})_{kl}=\delta_{ik}\delta_{jl}$ in a standard bases in
$Mat(2n,\mathbb{C})$. The evaluation homomorphism then defined by
the evaluation at a particular value of the variable $u$.

Then classical $L$-operator for $A_{2n-1}^{(2)}$ Toda chain has the
following form in this representation \be\label{Loperator}
L(p,x;u)\,=\,\sum_{i=1}^n
p_{x_i}\Big(E_{i,i}-E_{2n+1-i,2n+1-i}\Big)\,+
q_i(x)\Big(E_{i,i+1}-E_{2n-i,2n+1-i}\Big)-\\
\Big(E_{i+1,i}-E_{2n+1-i,2n-i}\Big)+
2q_0(x)\Big(E_{2n-1,1}-E_{2n,2}\Big)-\frac{u}{2}
\Big(E_{1,2n-1}-E_{2,2n}\Big),
\ee
where $q_i(x)=e^{(\alpha_i,x)}$ and $x=\sum_{i=1}^n x_ie_i$. 
 The space of  $L$-operators (\ref{Loperator})
provides a model for the phase space of the classical $A_{2n-1}^{(2)}$
Toda chain. The quantum $\cl$-operator is obtained by the standard
substitution $p_{x_j}=-i\hbar\pr_{x_j}$.
\be\label{qLoperator}
\cl(\pr_x,x;u)\,=\,\sum_{i=1}^n
-i\Big(E_{i,i}-E_{2n+1-i,2n+1-i}\Big)\,\frac{\pr}{\pr x_i}+
q_i(x)\Big(E_{i,i+1}-E_{2n-i,2n+1-i}\Big)-\\
\Big(E_{i+1,i}-E_{2n+1-i,2n-i}\Big)+
2q_0(x)\Big(E_{2n-1,1}-E_{2n,2}\Big)-\frac{u}{2}
\Big(E_{1,2n-1}-E_{2,2n}\Big).
\ee

 Classical Hamiltonian operators
of the  Toda chain are given by the coefficients of the
characteristic polynomial of $L$ \be \label{classhamgen}
\det(L(u)-\la)=u+\frac{1}{u}+\sum_{i=k}^N \,h_{2k}(p,x)\,
\la^{2n-2k} \ee For example the quadratic Hamiltonian generator is
given by  \be
H_2^{A_{2n-1}^{(2)}}(x)=\frac{1}{2}\sum_{i=1}^n\,p_{x_i}^2+e^{2x_1}+
\sum_{i=1}^{n-1}e^{x_{i+1}-x_i}\,+e^{-x_n-x_{n-1}} \ee
In \cite{GLO} we introduced  the following  integral  operator
defined by the  integral kernel
 \bqa\label{inttwA2odd}
Q_{A^{(2)}_{2n-1}}(x_i,\,z_i)&=& \exp (\frac{i}{\hbar}\mathcal{F}(x_i,z_i))=\\
&=&\exp\,\frac{i}{\hbar}\,\Big\{\, g_1e^{x_1+z_1}+ \sum_{i=1}^{n-1}\Big(
e^{x_i-z_i}+g_{i+1}e^{z_{i+1}-x_i}\Big)+
e^{x_n-z_n}+g_{n+1}e^{-x_{n}-z_n}\,\Big\},\nonumber\eqa
The action of the $Q$-operator on the functions is given by
\be
(Q_{A_{2n-1}^{(2)}} f)(x_i)=\int \bigwedge_{i=1}^n dz_i
 \,\, Q_{A_{2n-1}^{(2)}}(x_i,z_i)\,f(z_i)
\ee
Note that $\mathcal{F}$ is a generating function
of the canonical transformation corresponding to the integral
operator in the classical limit.

It was shown in \cite{GLO} that thus defined integral operator
 (\ref{inttwA2odd}) intertwines quadratic Hamiltonian operators
$H_2^{A_{2n-1}^{(2)}}$ for different coupling constants
\be\label{intertw}
H_2(x_i,\pr_{x_i},g_i)\,Q(x_i,z_i)=Q(x_i,z_i,g_i)H_2(z_i,\pr_{z_i},g_i')
\ee
where we use the quadratic Hamiltonian operators 
\be
H_2^{A_{2n-1}^{(2)}}(x)=\frac{1}{2}\sum_{i=1}^n\,p_{x_i}^2+2g_1e^{2x_1}+
\sum_{i=1}^{n-1}g_i \,e^{x_{i+1}-x_i}\,+g_ng_{n+1}e^{-x_n-x_{n-1}} \ee
and the relation holds for  $g_i=g'_{n+2-i}=1$,  $i\neq 1$ and 
$g_1=g'_{n+1}=2$ (see \cite{GLO} for details). In the following 
we will assume that these particular coupling constants are chosen. 
Also note that in (\ref{intertw})  and in the following similar
identities  we shall assume that 
 the Hamiltonian operator on l.h.s. acts on the right and
the Hamiltonian on r.h.s. acts on the left.

Let us show that the integral operator $Q$ commutes with all Hamiltonian operators
of $A_{2n-1}^{(2)}$ Toda chain. The prove is based on the following  
result.  
\begin{te}\label{TheoremAn} 
Let $\cl(x_i,\pr_{x_i},g_i,u)$ be a quantum $L$-operator given by
(\ref{qLoperator}) and $Q(x_i,z_i)$ be an integral kernel
(\ref{inttwA2odd}). Then the following relation holds
\be
R(x_i,z_i,g_i,g'_i,u)\cl(x_i,\pr_{x_i},g_i,u)
Q(x_i,z_i,u)=Q(x_i,z_i,u)\cl(z_i,\pr_{z_i},g'_i,u) R(x_i,z_i,g_i,g'_i,u)
\nonumber \ee
where matrix $R$ is given by
\be \label{Rmatrix}
R(x_i,z_i,g_i,g'_i,u)=\Big(e^{x_1+z_1}\,E_{n,1}+\frac{u}{2}E_{1,n}\Big)+
\sum_{i=0}^n\Big(e^{z_i-x_{i-1}}\,E_{n+1,i+1}-uE_{i,n+i}\Big)+\\
\Big(e^{-x_n-z_n}\,E_{2n,n+1}+\frac{u}{2}E_{n+1,2n}\Big)+
\sum_{i=1}^n\Big(e^{x_{n+1-i}-z_{n+1-i}}\,E_{i,n+i}-E_{n+i,i}\Big)\ee 
\end{te}
The proof is straightforward and reduces to the verification of the
matrix identity \be \label{classicalmult}
R(x_i,z_i,g_i,g'_i,u)L(x_i,p_{x_i},g_i,u)
=L(z_i,p_{z_j},g'_i,u)R(x_i,z_i,g_i,g'_i,u) 
\ee where  \bqa\label{cantrans} 
p_{x_1}=\frac{\partial \mathcal{F}}{\partial x_1}=
e^{z_1+x_1}+e^{x_1-z_1}-e^{z_2-x_1},\qquad
p_{x_i}=\frac{\partial \mathcal{F}}{\partial x_i}=e^{x_i-z_i}-e^{z_{i+1}-x_i}\,,\\
p_{z_i}=-\frac{\partial \mathcal{F}}{\partial z_i}=
e^{x_i-z_i}-e^{z_i-x_{i-1}},\qquad p_{z_n}=-\frac{\partial
\mathcal{F}}{\partial
z_n}=e^{-x_n-z_n}+e^{x_n-z_n}-e^{z_n-x_{n-1}}\eqa 

Note that the identity (\ref{classicalmult}) follows from the following
representations of  the pair of $L$-operators
\be\label{factone}
L(x_i,p_{x_i},g_i,u)=R\,\,R^*,\,\,\,\,\,
L(z_i,p_{z_i},g'_i,u)=R^*\,\,R, \ee
where $R$ is given by (\ref{Rmatrix})  and relations
(\ref{cantrans}) are implied. The involution $R\to R^*$ is the
involution that describes the twisted affine Lie algebra
$A_{2n-1}^{(2)}$ as a fixed point subalgebra of $A_{2n-1}^{(1)}$ 
(see \cite{DS} for the details). The equations (\ref{factone})
implies that $L$ is in the image of $A_{2n-1}^{(2)}$.    

We would like also note that the similar type of relations 
was used in the theory of Darboux maps (see \cite{AvM},\cite{V}).

\begin{ex}  In the case of $n=4$ adopt the following notations
$$a_i=e^{x_i-z_i},\qquad d_0=e^{x_1+z_1},\quad
d_i=e^{z_i-x_{i-1}},\quad d_n=e^{-x_n-z_n}\,.$$ In these terms we
get for $L(x_i,p_{x_i},g_i,u)$ with $g_i=1$
$$
\left(\begin{array}{cccccccc}
a_4-d_4 & a_4d_3 & 0 & 0 & 0 & 0 & -u/2 & 0\\
-1 & a_3-d_3 & a_3d_2 & 0 & 0 & 0 & 0 & -u/2\\
0 & -1 & a_2-d_2 & a_2d_1 & 0 & 0 & 0 & 0\\
0 & 0 & -1 & a_1+d_0-d_1 & 2a_1d_0 & 0 & 0 & 0\\
0 & 0 & 0 & -1 & -a_1-d_0+d_1 & 2a_2d_1 & 0 & 0\\
0 & 0 & 0 & 0 & -1 & -a_2+d_2 & a_3d_2 & 0\\
2u^{-1}d_3d_4 & 0 & 0 & 0 & 0 & -1 & -a_3+d_3 & a_4d_3\\
0 & 2u^{-1}d_3d_4 & 0 & 0 & 0 & 0 & -1 & -a_4+d_4
\end{array}\right).$$
Then we have
 \bqa R=\, \left(\begin{array}{cccccccc}
0 & 0 & 0 & -u/2 & e^{x_4-z_4}u & 0 & 0 & 0\\
0 & 0 & 0 & 0 & -u & e^{x_3-z_3}u & 0 & 0\\
0 & 0 & 0 & 0 & 0 & -u & e^{x_2-z_2}u & 0\\
2e^{z_1+x_1} & 0 & 0 & 0 & 0 & 0 & -u & e^{x_1-z_1}u\\
-1 & e^{z_2-x_1} & 0 & 0 & 0 & 0 & 0 & -u/2\\
0 & -1 & e^{z_3-x_2} & 0 & 0 & 0 & 0 & 0\\
0 & 0 & -1 & e^{z_4-x_3} & 0 & 0 & 0 & 0\\
0 & 0 & 0 & -1 & 2e^{-z_4-x_4} & 0 & 0 & 0
\end{array}\right)\eqa and
$$R^*\,=\, \left(\begin{array}{cccccccc}
0 & 0 & 0 & 1/2 & e^{x_1-z_1} & 0 & 0 & 0\\
0 & 0 & 0 & 0 & 1 & e^{x_2-z_2} & 0 & 0\\
0 & 0 & 0 & 0 & 0 & 1 & e^{x_3-z_3} & 0\\
2u^{-1}e^{-z_4-x_4} & 0 & 0 & 0 & 0 & 0 & 1 & e^{x_4-z_4}\\
u^{-1} & u^{-1}e^{z_4-x_3} & 0 & 0 & 0 & 0 & 0 & 1/2\\
0 & u^{-1} & u^{-1}e^{z_3-x_2} & 0 & 0 & 0 & 0 & 0\\
0 & 0 & u^{-1} & u^{-1}e^{z_2-x_1} & 0 & 0 & 0 & 0\\
0 & 0 & 0 & u^{-1} & 2u^{-1}e^{z_1+x_1} & 0 & 0 & 0
\end{array}\right)$$
\end{ex}
Taking into account (\ref{classicalmult}) we have
the following proposition. 
\begin{prop} The following identity holds
 \bqa\label{detident}
\det(L(x_i,p_{x_i},g_i,u)-\lambda)\,=\,
\det(L(z_i,p_{z_i},g'_i,u)-\lambda)\eqa
\end{prop}
The same reasoning goes for the quantum $L$-operators. 
Let $H_i^{A_{2n-1}^{(2)}}(x_i,\pr_{x_i},g_i)$ be  a full set of the
commuting quantum Hamiltonian operators for $A_{2n-1}^{(2)}$ Toda chain
given by the quantization of the classical Hamiltonian
operators (\ref{classhamgen}). The quantization uses the 
recursive representation for the characteristic polynomial  of the 
quantum $L$-operator and is based on the special properties of the matrix
differential operator (\ref{qLoperator}). 
 Then the  following intertwining relations holds
\be\label{qdetident}
H_i^{A_{2n-1}^{(2)}}(x_i,\pr_{x_i},g_i)Q(x_i,z_i)=
Q(x_i,z_i) H_i^{A_{2n-1}^{(2)}}(z_i,\pr_{z_i},g'_i) \nonumber
\ee
Let us stress that the  constructed integral  operators
intertwine the Hamiltonians of the $A_{2n-1}^{(2)}$ Toda chain 
with different coupling constants. It was explained in \cite{GLO} that
to get the Baxter $Q$-operator for $A_{2n-1}^{(2)}$ commuting with the
Hamiltonians one should apply the elementary intertwining twice 
 \bqa
Q^{A^{(2)}_{2n-1}}(x_1,\ldots,x_n;\,y_1,\ldots,y_n)=\int
\bigwedge_{i=1}^{n}dz_i\,
Q(x_1,\ldots,x_n;\,z_1,\ldots,z_{n})\times
\\ \nonumber
\times Q(y_1,\ldots,y_{n};\,z_1,\ldots,z_n).\eqa

\section{Recursive operators for $C_n$ and $D_n$}

In this section we prove the intertwining 
relations between recursive operators for 
the  full set of Hamiltonian operators for 
for $C_n$ and $D_n$ open Toda chains. 
These commutation relations can be also obtained 
from the analogous properties of Baxter $Q$-operator
 of $A_{2n-1}^{(2)}$ Toda chain in a certain limit.

We start with the description of $C_n$ and $D_n$ Toda chains.
The root system $C_n$ and $D_n$ and the Hamiltonian 
are defined in terms of the  orthonormal  basis $\{e_1,\ldots,e_{n}\}$ in
$\mathbb{R}^{n}$ as follows.
Simple roots of the root system $C_n$ are given by  \be \alpha_1=2e_1,\qquad
\alpha_{i+1}=e_{i+1}-e_i,\qquad 1\leq i\leq n-1\qquad
\alpha_{n+1}=-e_n-e_{n-1}, \ee
The corresponding Dynkin diagram is
$$
\begin{CD}\alpha_1 \Longrightarrow \alpha_2 @>>> \ldots @>>> \alpha_{n-1}
@>>> \alpha_n\\ @. @. @VVV @.\\ @. @.\alpha_{n+1}@.\end{CD}
$$
Simple roots of the root system $D_n$ are given by  \be \alpha_1=2e_1,\qquad
\alpha_{i+1}=e_{i+1}-e_i,\qquad 1\leq i\leq n-1\qquad
\alpha_{n+1}=-e_n-e_{n-1}, \ee
The corresponding Dynkin diagram is
$$
\begin{CD}\alpha_1 \Longrightarrow \alpha_2 @>>> \ldots @>>> \alpha_{n-1}
@>>> \alpha_n\\ @. @. @VVV @.\\ @. @.\alpha_{n+1}@.\end{CD}
$$
Classical Hamiltonian operators $h^{C_n}_i$ of $C_n$ Toda chain are given by the
coefficients of the characteristic polynomial
\be \label{charpol}
\det(L(x,p)-\la)=\sum_{i} h_i\la^{n-i}
\ee
where $L$ operator is given by
\be
L^{C_n}(x,p)=\sum_{i=1}^n\,p_i\Big(E_{i,i}-E_{2n+1-i,2n+1-i}\Big)+\\
q_i(x)\Big(E_{i,i+1}+E_{2n-i,2n+1-i}\Big)-
\Big(E_{i+1,i}-E_{2n+1-i,2n-i}\Big) \ee
Similarly classical Hamiltonian operators of $D_n$ Toda chain are
given by the coefficients of characteristic polynomial (\ref{charpol})
of the following $L$-operator
\be
 L^{D_n}(x,p)=
\sum_{i=1}^n\,p_i\Big(E_{i,i}-E_{2n+1-i,2n+1-i}\Big)+\\
\sum_{i=1}^{n-1}\,q_i(x)\Big(E_{i,i+1}+E_{2n-i,2n+1-i}\Big)-
\Big(E_{i+1,i}-E_{2n+1-i,2n-i}\Big)+\\
2q_n(x)\Big(E_{n-1,n+1}+E_{n,n+2}\Big)-(E_{n+1,n-1}+E_{n+2,n}) \ee

The corresponding quantum Hamiltonians are obtained by the standard
quantisation procedure. Thus for the quadratic Hamiltonian operators
we have
\be H_2^{C_n}(z_i)=-\frac{1}{2}\sum_{i=1}^n
\frac{\partial^2}{\partial z_i^2}\,+\sum_{i=1}^{n-1}
e^{z_{i+1}-z_i}\,+e^{-2z_n},\\
H_2^{D_n}(x_i)= -\frac{1}{2}\sum_{i=1}^n \frac{\partial^2}{\partial
x_i^2}\,+\sum_{i=1}^{n-1} e^{x_{i+1}-x_i}\,+e^{-z_n-z_{n-1}}\ee

In \cite{GLO} the integral operators intertwining
quadratic Toda chain Hamiltonian operators for different
elements of  the classical series $C_n$ and $D_n$ were proposed.
Integral operator defined by the kernel  \be\label{DnQCn}
Q^{D_{n}}_{C_{n}}(x_i;z_i)\,=\exp\, \frac{i}{\hbar}\mathcal{F}^{D_{n}}_{C_{n}}(x_i;z_i)
=\exp\,\frac{i}{\hbar}\Big\{\,\sum_{i=1}^{n-1}\Big(e^{x_i-z_i}+
e^{z_{i+1}-x_i}\Big)\,+\,e^{x_n-z_n}+e^{-x_n-z_n}\,\Big\},\ee
satisfies  the following intertwining relation
\be \label{HamintCn}
H_2^{C_n}(x_i,\pr_{x_i})\,Q_{D_{n}}^{C_{n}}(x_i,z_i)=
Q_{D_{n}}^{C_{n}}(x_i,z_i)H_2^{D_n}(z_i,\pr_{z_i}).
\ee
Similarly integral operator defined by the following kernel  \be \label{DnCnminus}
 Q_{D_n}^{\,\,\,C_{n-1}}(x_i;z_i)=\exp\,\frac{i}{\hbar}
 \mathcal{F}_{D_{n}}^{C_{n-1}}(x_i;z_i)=
\exp\,\frac{i}{\hbar}\Big\{\sum_{i=1}^{n-1}\Big(e^{z_i-x_i}+g_ie^{x_{i+1}-z_i}\Big)+
g_ne^{-x_n-z_{n-1}}\Big\},\ee
satisfies  the following intertwining relation
\be \label{HamintDn}
H_2^{D_n}(x_i,\pr_{x_i})\,Q_{D_{n}}^{C_{n-1}}(x_i,z_i)=
Q_{D_{n}}^{C_{n-1}}(x_i,z_i)H_2^{C_{n-1}}(z_i,\pr_{z_i})
\ee
In the classical approximation (\ref{HamintCn}) 
is   reduced to the identity 
\be
H_2^{C_n}(x_i,p_{x_i})\,=H_2^{D_n}(z_i,p_{z_i})
\ee
where the following expressions for the  momenta variables
are implied
 \be\label{CanTransf}
p_{x_i}=\frac{\partial \mathcal{F}^{D_{n}}_{C_{n}}} {\partial x_i}=
e^{x_i-z_i}-e^{z_{i+1}-x_i},\\ p_{x_n}=\frac{\partial
\mathcal{F}^{D_{n}}_{C_{n}}} {\partial
x_n}=e^{x_n-z_n}-e^{-z_{n}-x_n}\qquad 1\leq i< n,\ee \be
p_{z_1}=-\frac{\partial \mathcal{F}^{D_{n}}_{C_{n}}}{\partial z_1}=
e^{x_1-z_1},\qquad p_{z_i}=-\frac{\partial
\mathcal{F}^{D_{n}}_{C_{n}}}{\partial z_i}=
e^{x_i-z_i}-e^{z_i-x_{i-1}},\\ p_{z_n}=-\frac{\partial
\mathcal{F}^{D_{n}}_{C_{n}}}{\partial z_n}=
e^{x_n-z_n}+e^{-z_n-x_n}-e^{z_n-x_{n-1}},\qquad 1<i<n.\ee 
Similarly  (\ref{HamintDn}) is reduced in the classical approximation 
to   \be
H_2^{D_n}(x_i,p_{x_i})\,= H_2^{C_{n-1}}(z_i,p_{z_i}) \ee where the
following expressions for the  momenta variables are implied \be 
\label{CanTransfone}
 p_{x_1}=\frac{\partial \mathcal{F}_{D_{n}}^{C_{n-1}}}
{\partial x_1}=-e^{z_1-x_1},\qquad p_{x_i}=\frac{\partial
\mathcal{F}_{D_{n}}^{C_{n-1}}}{\partial x_i}=
e^{x_i-z_{i-1}}-e^{z_i-x_i},\\ p_{x_n}=\frac{\partial
\mathcal{F}_{D_{n}}^{C_{n-1}}}{\partial x_n}=
e^{x_n-z_{n-1}}-e^{-z_{n-1}-x_n},\qquad 1<i<n,\ee \be
p_{z_i}=-\frac{\partial\mathcal{F}_{D_{n}}^{C_{n-1}}}{\partial
z_i}=e^{x_{i+1}-z_i}-e^{z_i-x_i},\\ p_{z_{n-1}}=\frac{\partial
\mathcal{F}_{D_{n}}^{C_{n-1}}}{\partial
z_{n-1}}=e^{-z_{n-1}-x_n}+e^{x_n-z_{n-1}}-e^{z_{n-1}-x_{n-1}},\qquad
1\leq i<n-1.\ee  

Toda chains corresponding to the root systems $C_n$ and $D_n$ can be
obtained from $A_{2n-1}^{(2)}$ Toda chain using an appropriate
limiting procedure. To do this explicitly one should 
explicitly introduce the coupling constant in Toda chains. 
 It was shown in \cite{GLO} that the operator \bqa
Q_{A^{(2)}_{2n-1}}(x_i,\,z_i)&=&\frac{i}{\hbar} \exp \mathcal{F}(x_i,z_i)=\\
&=&\exp \,\frac{i}{\hbar}\Big\{\, g_1e^{x_1+z_1}+ \sum_{i=1}^{n-1}\Big(
e^{x_i-z_i}+g_{i+1}e^{z_{i+1}-x_i}\Big)+
e^{x_n-z_n}+g_{n+1}e^{-x_{n}-z_n}\,\Big\}\nonumber\eqa satisfies  the
following relation
\be H_2^{A_{2n-1}^{(2)}}(x_i)Q_{A^{(2)}_{2n-1}}(x_i,\,z_i)\,=\,
Q_{A^{(2)}_{2n-1}}(x_i,\,z_i)\widetilde{H}_2^{A_{2n-1}^{(2)}}(z_i)
\ee
for the quadratic Hamiltonians associated to the root system of type
$A_{2n-1}^{(2)}$:
\be H_2^{A_{2n-1}^{(2)}}(x_i)= -\frac{1}{2}\sum_{i=1}^n
\frac{\partial^2}{\partial x_i^2}\,+2g_1e^{2x_1}+\sum_{i=1}^{n-1}
g_{i+1}e^{x_{i+1}-x_i}\,+g_ng_{n+1}e^{-x_n-x_{n-1}}
\ee and
\be
\widetilde{H}_2^{A_{2n-1}^{(2)}}(z_i)= -\frac{1}{2}\sum_{i=1}^n
\frac{\partial^2}{\partial z_i^2}\,+ g_1g_2e^{z_1+z_2}+
\sum_{i=1}^{n-1} g_{i+1}e^{z_{i+1}-z_i}\,+2g_{n+1}e^{-2z_n}
\ee
Taking the limit $g_1\rightarrow0$ we obtain that  $Q$-operator
\be
Q^{D_{n}}_{C_{n}}(x_i;z_i)\,=\,\lim_{g_1\rightarrow0}
Q_{A^{(2)}_{2n-1}}(x_i,z_i)
\ee  intertwines
\be H_2^{C_n}(z_i)\,=\,\lim_{g_1\rightarrow0}
\widetilde{H}_2^{A_{2n-1}^{(2)}}(z_i)
\ee and
\be H_2^{D_n}(x_i)\,=\,\lim_{g_1\rightarrow0}
H_2^{A_{2n-1}^{(2)}}(x_i)\,.
\ee These Hamiltonians are naturally
associated with  root systems of type $C_n$ and $D_n$ respectively that
are isomorphic to certain root sub-systems in the affine root system
 of type $A_{2n-1}^{(2)}$. Thus this limit procedure
provides the $L$-operators associated to the root sub-systems of
types $C_n$ and $D_n$.

Note that the limiting  operator $Q^{D_{n}}_{C_{n}}(x_i;z_i)$
in general intertwines Hamiltonians of  Toda chains 
with different coupling constants except last ones. Thus the following
 operators :
\be H_2^{C_n}(z_i)=-\frac{1}{2}\sum_{i=1}^n
\frac{\partial^2}{\partial z_i^2}\,+\sum_{i=1}^{n-1}
g'_{i+1}e^{z_{i+1}-z_i}\,+2g_ne^{-2z_n}
\ee and
\be H_2^{D_n}(x_i)=-\frac{1}{2}\sum_{i=1}^n
\frac{\partial^2}{\partial x_i^2}\,+\sum_{i=1}^{n-1}
g_{i+1}e^{x_{i+1}-x_i}\,+g_{n-1}g_ne^{-z_n-z_{n-1}}
\ee 
are intertwined by the integral operator with the kernel 
 \bqa Q^{D_{n}}_{C_{n}}(x_i;z_i)\,=\,
\exp\Big\{\,\sum_{i=1}^{n-1}\Big(\gamma_ie^{x_i-z_i}+
\beta_{i}e^{z_{i+1}-x_i}\Big)\,+\,
e^{x_n-z_n}+g_ne^{-x_n-z_n}\,\Big\}\nonumber\eqa where 
\be \gamma_i=\prod_{k=i}^{n-1}\frac{g'_k}{g_k},\qquad
\beta_i=g_i\prod_{k=i+1}^{n-1}\frac{g_k}{g'_k}\qquad 1\leq i<n\,.
\ee
Taking the limit $g'_1\rightarrow0$ we obtain
\be H_2^{C_{n-1}}(z_i)\,=\,\lim_{g'_1\rightarrow0}H_2^{C_n}(z_i)
\ee
Thus the  intertwiner 
\be Q^{D_n}_{\,\,\,C_{n-1}}(x_i;z_i)\,=\,
\lim_{g'_1\rightarrow0}Q^{D_{n}}_{C_{n}}(x_i;z_i)
\ee satisfying
\be H_2^{D_n}(x_i)\,Q^{D_{n}}_{C_{n-1}}(x_i,z_i)=
Q^{D_{n}}_{C_{n-1}}(x_i,z_i)H_2^{C_{n-1}}(z_i)
\ee
Therefore all the recursive operators can be obtained 
from $A_{2n-1}^{(2)}$ Baxter $Q$-operators in different limits. 
This implies that the generalization of the intertwining relations
(\ref{HamintCn}) and (\ref{HamintDn})  to the
 full set of the Hamiltonians of $C_n$ and $D_n$ open Toda chains 
can be obtained  from (\ref{detident}), (\ref{qdetident})  
 using the limiting procedure described above. 
The following  commutation relations with $L$ operators
are the analog of the Theorem \ref{TheoremAn}. 
\begin{te} Let $\cl^{C_n}(x_i,\pr_{x_i},g_i,u)$ be a quantum $L$-operator
of $C_n$ Toda chain and  
\\ let   $\cl^{D_n}(x_i,\pr_{x_i},g_i,u)$ be a quantum $L$-operator
of $D_n$ Toda chain. Then the following intertwining relations for the
kernels (\ref{DnQCn}),(\ref{DnCnminus}) of the integral operators hold
\bqa \label{basicidentCD}
M(x_i,z_i,g_i,g'_i,u)\cl^{C_n}(x_i,\pr_{x_i},g_i,u)
Q^{D_{n}}_{C_{n}}
(x_i,z_i,u)=\\
=Q^{D_{n}}_{C_{n}}
(x_i,z_i)\cl^{D_n}(z_i,\pr_{z_i},g'_i,u) M(x_i,z_i,g_i,g'_i,u)
\nonumber \eqa
\bqa \label{basicidentDC}
N(x_i,z_i,g_i,g'_i,u)\cl^{D_n}(x_i,\pr_{x_i},g_i,u)
Q_{D_n}^{C_{n-1}}(x_i,z_i,u)=\\
= Q_{D_n}^{C_{n-1}}(x_i,z_i)\cl^{C_{n-1}}(z_i,\pr_{z_i},g'_i,u) N(x_i,z_i,g_i,g'_i,u)
\nonumber \eqa
where matrices  $M$ and $N$ are given
 \be
M(x_i,z_i,g_i,g'_i,u)=\sum_{i=1}^n\Big(E_{i,i}+e^{z_i-x_{i-1}}\,E_{i,i+1}\Big)+\\
\Big(e^{x_{n+1-i}-z_{n+1-i}}\,E_{2n+1-i,2n+1-i}+E_{2n+1-i,2n-i}\Big) 
+\Big(e^{-x_n-z_n}E_{n,n+1}-\frac{1}{2}E_{n+1,n}\Big)

\ee \be
N(x_i,z_i,g_i,g'_i,u)=\sum_{i=1}^n\Big(e^{x_i-z_i}\,E_{i,i}-E_{i+1,i}\Big)\\
\Big(E_{2n+1-i,2n+1-i}-e^{z_{n+1-i}-z_{n-i}}\,E_{2n-i,2n+1-i}\Big) 
+\Big(e^{-x_n-z_n}E_{n,n+1}+\frac{1}{2}E_{n+1,n}\Big)
\ee

\end{te}
Similar to the case of $A_{2n-1}^{(2)}$ the  commutation 
relations (\ref{basicidentCD}),(\ref{basicidentDC})
provide the quantum intertwining relation 
with $C_n$ and $D_n$ Hamiltonian operators. 
This leads to the proof of the  zero eigenvalue property for the full sets of the
Hamiltonian operators applied to the corresponding wave functions \cite{GLO}.

\section{Integral representation of generic $D_n$  eigenfunctions}

In \cite{GLO}    integral representations for
 zero eigenvalue wave functions of $D_n$  Toda chain quadratic Hamiltonian
operators were  given. In this  Section we generalize these representations
to the case of generic eigenvalues.  

Following \cite{GLO} we represent the wave functions of $D_n$ Toda
chain using the recursive integral operators. Let us define an 
integral recursive operator by the following kernel  
 \bqa\label{intrepnonzero}
Q_{D_{k+1}}^{\,\,\,\,D_k}(x_{k+1,i},\,x_{k,j},\,\la_{k+1})=\int\bigwedge_{i=1}^k\,dz_{k,i}
\exp\Big\{\i\la_{k+1}(\sum_{i=1}^{k+1} x_{k+1,i}+\sum_{i=1}^{k}x_{k,i}-
2\sum_{i=1}^kz_{k,i})\Big\}\,\times \nonumber \\ \times 
\left(e^{-x_{k+1,k+1}}+e^{-x_{k,k}}\right)^{2\i\,\la_{k+1}}
\times\\ \times  Q_{D_{k+1}}^{\,\,\,\,C_k}(x_{k+1,1},\ldots,x_{k+1,k+1};
z_{k,1},\ldots,z_{k,k})\, Q_{C_{k}}^{\,\,\,\,D_k}
(z_{k,1},\ldots,z_{k,k};x_{k,1},\ldots,x_{k,k}),\nonumber \eqa
where 
 \be Q_{D_{k+1}}^{\,\,\,C_{k}}(x_{k+1,i},\, z_{k,i})=
\exp\Big\{\sum_{i=1}^{k}\Big(e^{z_{k,i}-x_{k+1,i}}+e^{x_{k+1,i+1}-z_{k,i}}\Big)+
e^{-x_{k+1,k+1}-z_{k,k}}\Big\},\ee
\be
Q_{C_k}^{\,\,\,\,D_{k}}(z_{k,i},\,x_{k,i})=\exp\Big\{\,
\sum_{i=1}^{k}\Big(e^{x_{k,i}-z_{k,i}}+e^{z_{k,i+1}-x_{k,i}}\Big)+
e^{x_{k,k}-z_{k,k}}+e^{-x_{k,k}-z_{k,k}}\Big\}.\ee 
One can straightforwardly 
verify that  the following intertwining 
relations holds 
\be
H^{D_{k+1}}(x_{k+1,j})\,\,Q_{D_{k+1}}^{\,\,\,\,D_k}(x_{k+1,j},\,x_{k,j},\,\la_{k+1})= 
Q_{D_{k+1}}^{\,\,\,\,D_k}(x_{k+1,j},\,x_{k,j},\,\la_{k+1})\,
(H^{D_k}(x_{k,j})+\frac{1}{2}\la_{k+1}^2), \nonumber 
\ee
where  quadratic Hamiltonian operators are given by  
\be
 H^{D_k}(x_{k,j})=-
\frac{1}{2}\sum_{j=1}^k\frac{\partial^2}{\partial x_{k,j}^2}+
\sum_{i=1}^{k-1} e^{x_{k,j+1}-x_{k,j}}+e^{-x_{k,k}-x_{k,k-1}}. 
\ee
Consider  the  wave function for  $D_n$  Toda chain satisfying 
the eigenfunction equation 
\be 
 H^{D_n}(x)\,\,\,\Psi^{D_n}_{\la_1,\cdots ,\la_n} (x_1,\ldots,x_n)=
\frac{1}{2}(\sum_{i=1}^n\la_i^2)\,\,\,
\Psi^{D_n}_{\la_1,\cdots ,\la_n} (x_1,\ldots,x_n).
\ee
Then it can be represented in the following integral from
\be\label{WDN}
\Psi^{D_n}_{\la_1,\cdots ,\la_n} (x_1,\ldots,x_n)\,=\,
\int\bigwedge_{k=1}^{n-1}\bigwedge_{i=1}^kdx_{k,i}\,\,
e^{\i\la_1 x_{11}}\,\,
\prod_{k=1}^{n-1}Q_{D_{k+1}}^{\,\,\,\,D_k}(x_{k+1,j},\,x_{k,j},\,\la_{k+1}),\ee
where $x_i:=x_{n,i}$. 

Let us note that there is another form of  the integral representation 
for the kernel of the recursive operator (\ref{intrepnonzero})
 \be
Q_{D_{k+1}}^{\,\,\,\,D_k}(x_{k+1,j},\,x_{k,j},\,\la_{k+1})=
\frac{1}{\Gamma(2\i \,\la_{k+1})}\,
\int\bigwedge_{i=1}^k\,dz_{k,i}
\,\, d y_k\,\,\times \\ \times  
\exp\Big\{\i\la_{k+1}(\sum_{i=1}^{k+1} x_{k+1,i}+\sum_{i=1}^{k}x_{k,i}-
2\sum_{i=1}^kz_{k,i}+2y_{k+1})\Big\}\,\times \\ 
 \times  Q_{D_{k+1}}^{\,\,\,\,C_k}(x_{k+1,1},\ldots,x_{k+1,k+1};
z_{k,1},\ldots,z_{k,k};y_k)\, Q_{C_{k}}^{\,\,\,\,D_k}
(z_{k,1},\ldots,z_{k,k};x_{k,1},\ldots,x_{k,k};y_k).
\ee
where 
 \be Q_{D_{k+1}}^{\,\,\,C_{k}}(x_{k+1,i};\, z_{k,i},y_k)=
\exp\Big\{\sum_{i=1}^{k}\Big(e^{z_{k,i}-x_{k+1,i}}+e^{x_{k+1,i+1}-z_{k,i}}\Big)+
e^{-x_{k+1,k+1}-z_{k,k}}+e^{-x_{k+1,k+1}-y_k}\Big\},\nonumber \ee
\be
Q_{C_k}^{\,\,\,\,D_{k}}(z_{k,i};x_{k,i},y_k)=\exp\Big\{\,
\sum_{i=1}^{k}\Big(e^{x_{k,i}-z_{k,i}}+e^{z_{k,i+1}-x_{k,i}}\Big)+
e^{x_{k,k}-z_{k,k}}+e^{-x_{k,k}-z_{k,k}}+e^{-x_{k,k}-y_k}\Big\}.\nonumber
\ee
This integral representation  can be easily described in terms of the  
appropriate  Givental diagram. The detailed discussion of the relation
with the geometry of $D_n$  flag spaces and the results for other
classical  finite and affine Lie algebras  will be  
 presented in \cite{GLO1}.

\section{Appendix: Wave function for $D_2$}
Open Toda chain  for the  root system $D_2$ can be represented as a
couple of non-interacting $A_1$  Toda chains. This is a consequence of the 
isomorphism $D_2=A_1\oplus A_1$.  Quadratic Hamiltonian 
operator of $D_2$ Toda chain is given by 
 \be \label{d2quadratic}
H^{D_2}(x)=-\frac{1}{2}\sum\limits_{i=1}^{2}
\frac{\partial^2}{\partial x_i^2}+
e^{x_{22}-x_{21}}+e^{-x_{22}-x_{21}} \ee 
 Changing the variables \be \xi=x_{22}-x_{21},\,\,\,\,\,\, \eta =
x_{22}+x_{21}, \ee one transforms $H^{D_2}$  into the
sum of two $A_1$ Hamiltonians: \be \label{dtwo} H^{D_2}(x)=
2\,H^{A_1}(\xi)+2\,H^{A_1}(\eta):=
2 \Big(-\frac{1}{2} \frac{\partial^2}{\partial \xi^2}+
\frac{1}{2}e^{\xi}\Big)+2\Big(-\frac{1}{2} \frac{\partial^2}{\partial \eta^2}+
\frac{1}{2}e^{-\eta}\Big).\ee 
The additional quartic $D_2$ Hamiltonian commuting with
(\ref{d2quadratic}) is given by $H_4^{D_2}(x_1,x_2)=H^{A_1}(\xi)\,H^{A_1}(\eta)$.
Let $\Psi_{\la_1,\la_2}(x_{21},x_{22})$ be an eigenfunction 
of the quadratic $D_2$ Hamiltonian 
\be 
H^{D_2}(x)\,\,\Psi^{D_2}_{\la_1,\la_2}(x)=(\frac{1}{2}\la_1^2+\frac{1}{2}\la_2^2)
\,\,\Psi^{D_2}_{\la_1,\la_2}(x)
\ee
Taking into account (\ref{dtwo})  one can represent the eigenfunction 
as a product of  two eigenfunctions of $A_1$ Hamiltonians as
\be \label{factorisation}
\Psi_{\la_1,\la_2}(x_{21},x_{22})=\chi_{\frac{\la_2-\la_1}{2}}(-x_{22}-x_{21})\,
\chi_{\frac{\la_2+\la_1}{2}}(x_{22}-x_{21}), 
\ee
where $\chi_{\nu}(y)$ satisfies the equation 
\be \label{aone}
\Big(-\frac{1}{2} \frac{\partial^2}{\partial y^2}+
\frac{1}{2}e^{y}\Big)\chi_{\nu}(y)=\frac{1}{2}\nu^2\chi_{\nu}(y).
\ee
Solutions of (\ref{aone}) can expressed through Macdonald function as follows
\be 
\chi_{\nu}(y)=K_{2\i\nu}(2e^y)=-\frac{e^{2\pi \nu}
}{2}\int\limits_{-\infty+\i\pi}^{\infty+\i\pi}
e^{2\i\nu(x+y)}\exp\Big\{e^{-x}+e^{x+2y}\Big\}d x.
\ee 

Below  we explicitly demonstrate the factorisation 
(\ref{factorisation}) using the integral representation 
for $D_2$ Toda chain wave function described in the previous Section. 

The eigenfunction of $H^{D_2}$ has the  following integral
representation (see (\ref{WDN}) for $n=2$))
\be\label{psi} \Psi^{D_2}_{\la_1,\la_2}(x_{21},x_{22})=\int\int d
x_{11}d z_{11} dy \,\,\exp\Big\{
\i\la_2(x_{22}+x_{21})-\i\,2\,\la_2(z_{11}-y)+\i(\la_2+\la_1)x_{11}\Big\}\times
\\ \times
\exp\Big\{e^{z_{11}-x_{21}}+e^{x_{22}-z_{11}}+e^{-x_{22}-z_{11}}
+e^{x_{11}-z_{11}}+e^{-x_{11}-z_{11}}+
e^{-x_{22}-y}+e^{-x_{11}-y}\Big\}.
 \ee
One can integrate explicitly over $y$ to get the following
representation 
\be
\label{psinoy} \Psi^{D_2}_{\la_1,\la_2}(x_{21},x_{22})=\Gamma(2\i \la_2)\,\int\int d
x_{11}d z_{11} \,\,\left(e^{-x_{22}}+e^{-x_{11}}\right)^{\i
  \,2\,\la_2}\times \\ \times \exp\Big\{
\i\la_2(x_{22}+x_{21})-\i \,2\,\la_2\,z_{11}+\i(\la_2+\la_1)x_{11}\Big\}\times
\\ \times
\exp\Big\{e^{z_{11}-x_{21}}+e^{x_{22}-z_{11}}+e^{-x_{22}-z_{11}}
+e^{x_{11}-z_{11}}+e^{-x_{11}-z_{11}}\Big\}.
 \ee
 Let us introduce new variables 
\be
\xi=x_{22}-x_{21},\,\,\,\,\,\, \eta =
x_{22}+x_{21},\,\,\,\, z'_{11}=z_{11}-\frac{x_{22}+x_{21}}{2}, 
\,\,\,\,x'_{11}=x_{11}-\frac{x_{22}+x_{21}}{2} .\ee
Then \be \Psi^{D_2}_{\la_1,\la_2}(x_{21},x_{22})=\Gamma(2\i \la_2)\,\int\int d x_{11}' d
z_{11}'\,\,\left(e^{-x_{22}}+e^{-x'_{11}+\frac{\eta}{2}}\right)^{\i
  \,2\,\la_2}\times \\ \times \exp\Big\{
-\i\frac{(\la_2+\la_1)}{2}\eta-\i \,2\,\la_2\,z'_{11}+\i(\la_2+\la_1)x'_{11}\Big\}\times
\\ \times
\exp\Big\{
e^{z_{11}'+\frac{\xi}{2}}+e^{\frac{\xi}{2}-z_{11}'}+e^{x_{11}'-z_{11}'}+
e^{-\frac{\xi}{2}-\eta-z_{11}'}+e^{-x_{11}'-z_{11}'-\eta}\Big\}.
 \ee
Define  new variable as follows \be\label{anz}
e^{\tau}=e^{-z_{11}'-\eta}\Big(e^{-\frac{\xi}{2}}+e^{-x_{11}'}\Big),
 \ee
\be
e^{t}=e^{x'_{11}-z_{11}'-\frac{\xi}{2}-\eta}
\Big(e^{-\frac{\xi}{2}}+e^{-x_{11}'}\Big) .
\ee
As a result we  obtain the following integral representation 
 \be \Psi^{D_2}_{\la_1,\la_2}(x_{21},x_{22})=\Gamma(2\i \la_2)\,
\int\int d\tau d t \,\,
\exp\Big\{
-\i\frac{(\la_2-\la_1)}{2}\eta+\i \,(\la_2-\la_1)+\\+
\i\frac{(\la_2+\la_1)}{2}\xi + \i(\la_2+\la_1)t\Big\}
\exp \Big\{e^{-\tau}+e^{\tau-\eta}+e^{-t}
+e^{t+\xi}\Big\}=\\=4e^{-2\pi \la_2}\,\Gamma(2\i \la_2)\,\,
K_{\i\frac{\la_2+\la_1}{2}}(2e^{\frac{\xi}{2}})
K_{\i\frac{\la_2-\la_1}{2}}(2 e^{-\frac{\eta}{2}}).
\ee 
Thus we have demonstrated that the proposed integral representation 
for $D_n$  Toda chain wave functions  is compatible with the factorization 
of $D_2$ Toda chain.

\end{document}